\documentclass{amsart}

\usepackage{amsthm}
\usepackage[leqno]{amsmath}
\usepackage{latexsym,amsfonts,amssymb}
\usepackage{hyperref}

\newcommand{\numberseries}{\mdseries}   

\newlength{\thmtopspace}                
\newlength{\thmbotspace}                
\newlength{\thmheadspace}               
\newlength{\thmindent}                  

\newtheoremstyle{bfupright head,slanted body}
                {\thmtopspace}{\thmbotspace}
                {\slshape}{\thmindent}{\bfseries}{.}{\thmheadspace}
                {{\numberseries \thmnumber{(#2) }}\thmnote{#3}}

\newtheoremstyle{bfupright head,upright body}
                {\thmtopspace}{\thmbotspace}
                {\upshape}{\thmindent}{\bfseries}{.}{\thmheadspace}
                {{\numberseries \thmnumber{(#2) }}\thmnote{#3}}

\newtheoremstyle{fixed bf head,slanted body}
                {\thmtopspace}{\thmbotspace}{\slshape}
                {\thmindent}{\bfseries}{.}{\thmheadspace}
                {{\numberseries \thmnumber{(#2) }}\thmname{#1}\thmnote{ (#3)}}

\newtheoremstyle{fixed bf head,upright body}
                {\thmtopspace}{\thmbotspace}{\upshape}
                {\thmindent}{\bfseries}{.}{\thmheadspace}
                {{\numberseries \thmnumber{(#2) }}\thmname{#1}\thmnote{ (#3)}}

\newtheoremstyle{numbered paragraph}
                {\thmtopspace}{\thmbotspace}{\upshape}
                {\thmindent}{\upshape}{}{0pt}
                {{\numberseries \thmnumber{(#2) }}}

\newtheoremstyle{unnumbered paragraph}
                {\thmtopspace}{\thmbotspace}{\upshape}
                {\parindent}{\upshape}{}{0pt}

\theoremstyle{bfupright head,slanted body}
\newtheorem{res}{}[section]             \newtheorem*{res*}{}

\theoremstyle{bfupright head,upright body}
\newtheorem{bfhpg}[res]{}               \newtheorem*{bfhpg*}{}

\theoremstyle{fixed bf head,slanted body}
\newtheorem{thm}[res]{Theorem}          \newtheorem*{thm*}{Theorem}
\newtheorem{prp}[res]{Proposition}      \newtheorem*{prp*}{Proposition}
\newtheorem{cor}[res]{Corollary}        \newtheorem*{cor*}{Corollary}
\newtheorem{lem}[res]{Lemma}            \newtheorem*{lem*}{Lemma}

\theoremstyle{fixed bf head,upright body}
       \newtheorem*{dfn*}{Definition}
      \newtheorem*{obs*}{Observation}
\newtheorem{rmk}[res]{Remark}           \newtheorem*{rmk*}{Remark}
\newtheorem{exa}[res]{Example}          \newtheorem*{exa*}{Example}
            \newtheorem{stp*}{Setup}

\theoremstyle{numbered paragraph}
\newtheorem{ipg}[res]{}

\newlength{\thmlistleft}        
\newlength{\thmlistright}       
\newlength{\thmlistpartopsep}   
\newlength{\thmlisttopsep}      
\newlength{\thmlistparsep}      
\newlength{\thmlistitemsep}     

\newcounter{eqc} 
\newenvironment{eqc}{\begin{list}{\upshape (\textit{\roman{eqc}})}%
    {\usecounter{eqc}%
      \setlength{\leftmargin}{\thmlistleft}%
      \setlength{\labelwidth}{\thmlistleft}%
      \setlength{\rightmargin}{\thmlistright}%
      \setlength{\partopsep}{\thmlistpartopsep}%
      \setlength{\topsep}{\thmlisttopsep}%
      \setlength{\parsep}{\thmlistparsep}%
      \setlength{\itemsep}{\thmlistitemsep}}}%
  {\end{list}}%

\newcommand{\eqclbl}[1]{{\upshape(\textit{#1})}}

\newcounter{prt}
\newenvironment{prt}{\begin{list}{\upshape (\alph{prt})}%
    {\usecounter{prt}%
      \setlength{\leftmargin}{\thmlistleft}%
      \setlength{\labelwidth}{\thmlistleft}%
      \setlength{\rightmargin}{\thmlistright}%
      \setlength{\partopsep}{\thmlistpartopsep}%
      \setlength{\topsep}{\thmlisttopsep}%
      \setlength{\parsep}{\thmlistparsep}%
      \setlength{\itemsep}{\thmlistitemsep}}}%
  {\end{list}}%

\newcommand{\prtlbl}[1]{{\upshape(#1)}}

\newenvironment{itemlist}{\nopagebreak \begin{list}{$\bullet$}%
    {\setlength{\leftmargin}{\thmlistleft}%
      \setlength{\labelwidth}{\thmlistleft}%
      \setlength{\rightmargin}{\thmlistright}%
      \setlength{\partopsep}{\thmlistpartopsep}%
      \setlength{\topsep}{\thmlisttopsep}%
      \setlength{\parsep}{\thmlistparsep}%
      \setlength{\itemsep}{\thmlistitemsep}}}%
  {\end{list}}%

\newenvironment{prf*}[1][Proof]{%
  \begin{proof}[\bf #1]
    \setcounter{equation}{0}
    \renewcommand{\theequation}{\arabic{equation}}}
  {\end{proof}
}

\newcommand{\pgref}[1]{(\ref{#1})}
\newcommand{\thmref}[2][Theorem~]{#1\pgref{thm:#2}}
\newcommand{\corref}[2][Corollary~]{#1\pgref{cor:#2}}
\newcommand{\prpref}[2][Proposition~]{#1\pgref{prp:#2}}
\newcommand{\lemref}[2][Lemma~]{#1\pgref{lem:#2}}

\newcommand{\rmkref}[2][Remark~]{#1\pgref{rmk:#2}}
\newcommand{\secref}[2][Section~]{#1\ref{sec:#2}}
\newcommand{\partpgref}[2]{(\ref{#1})\prtlbl{#2}}

\newcommand{\partprpref}[3][Proposition~]{#1\partpgref{prp:#2}{#3}}
\renewcommand{\eqref}[1]{\pgref{eq:#1}}

\newcommand{\thmcite}[2][?]{\cite[thm.~#1]{#2}}
\newcommand{\corcite}[2][?]{\cite[cor.~#1]{#2}}
\newcommand{\prpcite}[2][?]{\cite[prop.~#1]{#2}}
\newcommand{\lemcite}[2][?]{\cite[lem.~#1]{#2}}
\newcommand{\dfncite}[2][?]{\cite[def.~#1]{#2}}
\newcommand{\seccite}[2][?]{\cite[sec.~#1]{#2}}

\newcommand{\set}[2][\mspace{1mu}]{\{#1 #2 #1\}}
\newcommand{\setof}[3][\mspace{2mu}]{\{#1#2 \mid #3#1\}}
\newcommand{\ZZ}{\mathbb{Z}}

\newcommand{\iinZ}{{i\in\ZZ}}
\newcommand{\jinZ}{{j\in\ZZ}}
\newcommand{\qtext}[1]{\quad\text{#1}\quad}
\newcommand{\qqtext}[1]{\qquad\text{#1}\qquad}
\newcommand{\qand}{\qtext{and}}
\newcommand{\qqand}{\qqtext{and}}

\DeclareMathOperator*{\dcoprod}{\textstyle\coprod}
\DeclareMathOperator*{\dprod}{\textstyle\prod}
\DeclareMathOperator*{\dfinprod}{\textstyle\bigoplus}

\newcommand{\m}{\mathfrak{m}}
\newcommand{\p}{\mathfrak{p}}

\newcommand{\is}{\cong}
\newcommand{\qis}{\simeq}
\renewcommand{\le}{\leqslant}
\renewcommand{\ge}{\geqslant}
\newcommand{\onto}{\twoheadrightarrow}

\newcommand{\lra}{\longrightarrow}
\newcommand{\xla}[2][]{\xleftarrow[#1]{\;#2\;}}
\newcommand{\xra}[2][]{\xrightarrow[#1]{\;#2\;}}
\newcommand{\qra}{\xra{\;\qis\;}}

\newcommand{\Rhat}{\widehat{R}}
\newcommand{\mapdef}[4][\rightarrow]{\nobreak{#2\colon #3 #1 #4}}
\newcommand{\qisdef}[4][\xra{\qis}]{\nobreak{#2\colon #3 #1 #4}}
\newcommand{\Ker}[1]{\nobreak{\operatorname{Ker}#1}}
\newcommand{\Coker}[1]{\nobreak{\operatorname{Coker}#1}}
\newcommand{\Cone}[1]{\nobreak{\operatorname{Cone}#1}}
\newcommand{\tev}[1]{\omega_{#1}}
\newcommand{\sign}[1]{(-1)^{#1}}
\newcommand{\dif}[2][]{{\partial}^{#2}_{#1}}
\newcommand{\Co}[2][]{\operatorname{C}_{#1}(#2)}
\renewcommand{\H}[2][]{\operatorname{H}_{#1}(#2)}
\renewcommand{\Im}[1]{\operatorname{Im}{#1}}
\newcommand{\Shift}[2][]{\mathsf{\Sigma}^{#1}{#2}}
\newcommand{\dptR}{\operatorname{depth}R}
\newcommand{\dimR}{\operatorname{dim}R}
\newcommand{\edim}[1]{\operatorname{edim}#1}
\newcommand{\codim}[1]{\operatorname{codim}#1}
\newcommand{\codepth}[1]{\operatorname{codepth}#1}
\newcommand{\wdt}[2][R]{\operatorname{width}_{#1}#2}
\newcommand{\dpt}[2][R]{\operatorname{depth}_{#1}#2}
\newcommand{\hgt}[2][R]{\operatorname{ht}_{#1}#2}
\newcommand{\id}[2][R]{\operatorname{id}_{#1}#2}
\newcommand{\pd}[2][R]{\operatorname{pd}_{#1}#2}
\newcommand{\CIdim}[2][R]{\operatorname{CI-dim}_{#1}#2}
\newcommand{\Gdim}[2][R]{\operatorname{G-dim}_{#1}#2}
\newcommand{\Gpd}[2][R]{\operatorname{Gpd}_{#1}#2}
\newcommand{\Hom}[3][R]{\operatorname{Hom}_{#1}(#2,#3)}
\newcommand{\RHom}[3][R]{\operatorname{\mathbf{R}Hom}_{#1}(#2,#3)}
\newcommand{\Ext}[4][R]{\operatorname{Ext}_{#1}^{#2}(#3,#4)}
\newcommand{\tp}[3][R]{\nobreak{#2\otimes_{#1}#3}}
\newcommand{\tpP}[3][R]{\left(\tp[#1]{#2}{#3}\right)}
\newcommand{\tpp}[3][R]{(\tp[#1]{#2}{#3})}
\newcommand{\Ltp}[3][R]{\nobreak{#2\otimes_{#1}^{\mathbf{L}}#3}}
\newcommand{\Ltpp}[3][R]{(\Ltp[#1]{#2}{#3})}
\newcommand{\Tor}[4][R]{\operatorname{Tor}^{#1}_{#2}(#3,#4)}

\numberwithin{equation}{res}

\newcommand{\btp}[3][R]{\nobreak{#2\otimes^\Join_{#1}#3}}

\newcommand{\Ttor}[4][R]{\smash{\operatorname{\widehat{Tor}}%
  }_{#2}^{#1^{\phantom{|\mspace{-6mu}}
    }}(#3,#4)}
\newcommand{\Text}[4][R]{\smash{\operatorname{\widehat{Ext}}%
  }_{#1}^{#2^{\phantom{|}\mspace{-6mu}}}(#3,#4)}

\newcommand{\Koszul}[2][R]{\operatorname{K}^{#1}(#2)}
\def\urltilda{\kern -.15em\lower .7ex\hbox{\~{}}\kern .04em}
 \newcommand{\Eterm}[3]{\operatorname{{}^{#1}E}^2_{#2,#3}}
 \newcommand{\Einfterm}[3]{\operatorname{{}^{#1}E}^\infty_{#2,#3}}

  {\end{list}}%

\begin{document}

\title[Vanishing of Tate homology and depth formulas over local
rings]{Vanishing of Tate homology\\ and depth formulas over local
  rings}

\author[L.\,W. Christensen]{Lars Winther Christensen}

\address{Department of Math.\ and Stat., Texas Tech University,
  Lubbock, TX 79409, U.S.A.}

\email{lars.w.christensen@ttu.edu}

\urladdr{http://www.math.ttu.edu/\urltilda lchriste}

\author[D.\,A. Jorgensen]{David A. Jorgensen}

\address{Department of Mathematics, University of Texas, Arlington,
  TX~76019, U.S.A.}

\email{djorgens@uta.edu}

\urladdr{http://dreadnought.uta.edu/\urltilda dave}

\thanks{This research was partly supported by NSA grants
  H98230-11-0214 (L.W.C.) and H98230-10-0197 (D.A.J.). Part of the
  work was done during the authors' visits to the University of
  Bielefeld (D.A.J) and the University of Paderborn (L.W.C.); the
  hospitality of both institutions is acknowledged with gratitude.}

\date{13 December 2013}

\dedicatory{To Hans-Bj{\o}rn Foxby on his 65$^{th}$ birthday}



\begin{abstract}
  Auslander's depth formula for pairs of Tor-independent modules over
  a regular local ring, \mbox{$\dpt[]{\tpp{M}{N}} = \dpt[]{M} +
    \dpt[]{N} - \dptR$}, has been generalized in several directions;
  most significantly it has been shown to hold for pairs of
  Tor-independent modules over complete intersection rings.

  In this paper we establish a depth formula that holds for every pair
  of Tate Tor-independent modules over a Gorenstein local ring.  It
  subsumes previous generalizations of Auslander's formula and yields
  new results on vanishing of cohomology over certain Gorenstein
  rings.
\end{abstract}

\maketitle
\thispagestyle{empty}

\section*{Introduction}
\label{sec:inro}

\noindent
To infer properties, qualitative or quantitative, of a tensor product
from properties of its factors is a delicate task. For finitely
generated modules $M$ and $N$ over a commutative noetherian local ring
$R$, Auslander \cite{MAs61} proved that the depth of the tensor
product is given by the formula
\begin{equation}
  \label{eq:df}\tag{A}
  \dpt{\tpp{M}{N}} = \dpt{M} + \dpt{N} - \dptR,
\end{equation}
provided that the projective dimension of $M$ is finite and the two
modules are Tor-independent, that is, the homology modules
$\Tor{i}{M}{N}$ vanish for $i \ge 1$. In particular, the equality
\eqref{df} holds for every pair of finitely generated Tor-independent
modules over a regular local ring.

The condition of finite projective dimension was first relaxed by
Huneke and Wiegand~\cite{CHnRWg94}, who established the validity of
\eqref{df} for pairs of finitely generated Tor-independent modules
over complete intersection local rings. Later, Araya and
Yoshino~\cite{TArYYs98} and Iyengar~\cite{SIn99} showed that \eqref{df}
holds for Tor-independent modules $M$ and $N$, provided that $M$ has
finite complete intersection dimension.

In a different direction, Foxby \cite{HBF80} relaxed the condition of
Tor-independence as follows. Let $M$ and $N$ be modules over a
commutative noetherian local ring $R$ and let $P$ be a projective
resolution of $M$. If $M$ has finite projective dimension, then there
is an equality,
\begin{equation}
  \label{eq:pdf}\tag{B}
  \dpt{\tpp{P}{N}} = \dpt{M} + \dpt{N} - \dptR.
\end{equation}
Implicit in this formula is an extension of the invariant ``depth'' to
complexes of modules; we recall it in \pgref{finite depth}. The
homology of the complex $\tp{P}{N}$ is $\Tor{*}{M}{N}$, and if $M$ and
$N$ are Tor-independent, then \eqref{pdf} reduces to Auslander's
formula \eqref{df}.
\begin{center}
  $\ast \ast \ast$
\end{center}
The notion of Tate homology for modules over group algebras has a
natural extension to modules over Gorenstein rings and, more
generally, to modules of finite Gorenstein projective dimension over
any ring. This theory was recently treated by Iacob~\cite{AIc07}. We
recall the basics in \pgref{tate}; a broader discussion is given
in~\seccite[2]{LWCDAJa}.

Let $R$ be a commutative noetherian local ring. The central result in
this paper, \thmref{dpt}, establishes vanishing of Tate homology
$\Ttor{*}{M}{N}$ as a sufficient condition for the equality
\eqref{pdf} to hold for a pair of $R$-modules $(M,N)$, where $M$ has
finite Gorenstein projective dimension. The Tate homology for such a
pair vanishes if $M$ has finite projective dimension over $R$, and if
$R$ is complete intersection, then Tor-independence implies vanishing
of Tate homology. Thus the main theorem subsumes all of the
aforementioned generalizations of Auslander's depth formula. In fact,
it goes further and subsumes several other generalizations obtained
over the half-century that has passed since \cite{MAs61} appeared.

What is more significant, though, is that \thmref{dpt} applies to
modules over all Gorenstein rings. The works of Iyengar and of Huneke
and Wiegand established the depth formula in the realm of complete
intersection rings. The depth is a cohomological invariant, and the
cohomological behavior of modules over Gorenstein rings can stray
dramatically from that of modules over complete intersection
rings. Therefore, it came as a surprise to us that vanishing of Tate
homology is sufficient to tame~depth in this vastly wider context.

In \secref[Sections~]{ci} and \secref[]{formula} we explore the
consequences of the main theorem for modules over different classes of
Gorenstein rings. For so-called AB rings, a notion coined by Huneke
and Jorgensen \cite{CHnDAJ03}, we obtain in \thmref{abbound} a precise
bound on vanishing of cohomology of finitely generated modules. For
modules over complete intersection rings, we obtain a derived depth
formula for modules that satisfy an effectively verifiable condition
on vanishing of Tate homology. That is, if $\Ttor{i}{M}{N}$ vanishes
for a finite number of consecutive indices---a number that only
depends on $R$---then (B) holds; this is \thmref{ci}.

In the final section we prove a statement, dual to \thmref{dpt}, for
the width invariant and use it to establish a bound on vanishing of
cohomology of finitely generated modules with vanishing Tate
cohomology $\Text{*}{M}{N}$.

\section{Depth of complexes}
\label{sec:complexes}

\noindent
In this paper, $R$-complexes---that is, complexes of $R$-modules---are
graded homologically.  A complex
\begin{equation*}
  M:\quad \cdots \lra M_{i+1} \xra{\dif[i+1]{M}} M_i \xra{\dif[i]{M}}
  M_{i-1} \lra \cdots
\end{equation*}
is called \emph{acyclic} if the homology complex $\H{M}$ is the
zero-complex. We use the notation $\Co[i]{M}$ for the cokernel of the
differential $\dif[i+1]{M}$. For $n\in\ZZ$ the \emph{n-fold shift} of
$M$ is the complex $\Shift[n]{M}$ given by $(\Shift[n]{M})_i =
M_{i-n}$ and $\dif[i]{\Shift[n]{M}} = (-1)^n\dif[i-n]{M}$.

The notation $\sup{M}$ and $\inf{M}$ is used for the supremum and
infimum of the set $\setof{i \in \ZZ}{M_i \ne 0}$, with the
conventions $\sup{\emptyset} = -\infty$ and $\inf{\emptyset} =
\infty$. A complex $M$ is called \emph{bounded above} if $\sup{M}$ is
finite, it is called \emph{bounded below} if $\inf{M}$ is finite, and
it is called \emph{bounded} if it is bounded above and below.

For some of our proofs, we shall need the following variation on
\lemcite[4.4.F]{LLAHBF91}.

\begin{lem}
  \label{lem:tev}
  Let $F$ and $M$ be $R$-complexes, and let $P$ be a complex of
  finitely generated $R$-modules. If one of the following conditions
  holds
  \begin{prt}
  \item $F$ and $M$ are bounded above, and $P$ is bounded below, or
  \item $M$ and $P$ are bounded
  \end{prt}
  and $P$ is a complex of projective $R$-modules or $F$ is a complex
  of flat $R$-modules, then there is a natural isomorphism of
  $R$-complexes
  \begin{equation*}
    \tp{\Hom{P}{M}}{F} \is \Hom{P}{\tp{M}{F}}.
  \end{equation*}
\end{lem}

\begin{prf*}
  For $R$-modules $P$, $M$, and $F$ the tensor evaluation map
  \begin{equation*}
    \mapdef{\tev{PMF}}{\tp{\Hom{P}{M}}{F}}{\Hom{P}{\tp{M}{F}}}
  \end{equation*}
  given by
  \begin{equation*}
    \tev{PMF}(\psi \otimes f)(p) = \psi(p)\otimes f,
  \end{equation*}
  is a homomorphism. It is an isomorphism if $P$ is finitely generated
  and projective, and also if $P$ is finitely generated and $F$ is
  flat; see \lemcite[4.4.F]{LLAHBF91}.

  Under either assumption, (a) or (b), the complex $M$ is bounded
  above and $P$ is bounded below. Assume, therefore, without loss of
  generality that one has $M_u=0$ for all $u>0$ and $P_u=0$ for all
  $u<0$. For every $n \in \ZZ$ one then has
  \begin{align*}
    \tpP{\Hom{P}{M}}{F}_n
    &= \dcoprod_\iinZ \tp{\Hom{P}{M}_i}{F_{n-i}} \\
    &= \dcoprod_{\iinZ} \tp{(\dprod_\jinZ
      \Hom{P_j}{M_{j+i}})}{F_{n-i}}\\
    &= \dcoprod_{i \le 0} \tp{(\dfinprod_{j=0}^{-i}
      \Hom{P_j}{M_{j+i}})}{F_{n-i}}\\[-.5\baselineskip]
    \intertext{and} \Hom{P}{\tp{M}{F}}_n
    &= \dprod_\jinZ \Hom{P_j}{\tpP{M}{F}_{j+n}} \\
    &= \dprod_{j \ge 0} \Hom{P_j}{\dcoprod_{h \in \ZZ}
      \tp{M_{h}}{F_{j+n-h}}} \\
    &= \dprod_{j \ge 0} \Hom{P_j}{\dcoprod_{i \le -j}
      \tp{M_{j+i}}{F_{n-i}}}.
  \end{align*}

  If (a) holds, then one can assume that $F_u$ is zero for all $u>0$,
  whence
  \begin{align*}
    \tpP{\Hom{P}{M}}{F}_n 
    &= \dfinprod_{i=n}^0
    \dfinprod_{j=0}^{-i}\tp{\Hom{P_j}{M_{j+i}}}{F_{n-i}}\\[-.5\baselineskip]
    \intertext{and} 
    \Hom{P}{\tp{M}{F}}_n 
    &= \dfinprod_{j=0}^{-n} \Hom{P_j}{\dfinprod_{i = n}^{-j}
      \tp{M_{j+i}}{F_{n-i}}}\\
    &= \dfinprod_{i=n}^0 \dfinprod_{j=0}^{-i}
    \Hom{P_j}{\tp{M_{j+i}}{F_{n-i}}}.
  \end{align*}
  The map from $\Hom{P}{\tp{M}{F}}$ to $\tp{\Hom{P}{M}}{F}$ with
  degree $n$ component $\dfinprod_{i=n}^0
  \dfinprod_{j=0}^{-i}\sign{j(n-i)}\tev{P_jM_{j+i}F_{n-i}}$ is a
  morphism of complexes; this is elementary to verify. It follows that
  it is an isomorphism if the modules in $P$ are projective or the
  modules in $F$ are~flat.

  If (b) holds, then a similar argument applies.
\end{prf*}

\begin{bfhpg}[Resolutions]
  \label{res}
  A morphism of $R$-complexes that induces an isomorphism in homology
  is called a \emph{quasi-isomorphism} and indicated by the symbol
  `$\qis$'.

  An $R$-complex $P$ is called \emph{semi-projective} if each module
  $P_i$ is projective, and the functor $\Hom{P}{-}$ preserves
  quasi-isomorphisms. Every bounded below complex of projective
  $R$-modules is semi-projective. Similarly, an $R$-complex $I$ is
  called \emph{semi-injective} if each module $I_i$ is injective, and
  the functor $\Hom{-}{I}$ preserves quasi-isomorphisms. Every bounded
  above complex of injective $R$-modules is semi-injective. The
  following facts are proved in \cite{dga}.
  \begin{prt}
  \item[(P)] Every $R$-complex $M$ has a semi-projective resolution.
    That is, there is a quasi-isomorphism $\mapdef{\pi}{P}{M}$, where
    $P$ is a semi-projective complex with $P_i=0$ for all $i <
    \inf{M}$. Moreover, if $\H{M}$ is bounded below, then $M$ has a
    semi-projective resolution $P' \qra M$ with $P'_i=0$ for all $i <
    \inf{\H{M}}$.

  \item[(I)] Every $R$-complex $M$ has a semi-injective resolution.
    That is, there is a quasi-isomorphism $\mapdef{\iota}{M}{I}$,
    where $I$ is semi-injective with $I_i=0$ for all $i >
    \sup{M}$. Moreover, if $\H{M}$ is bounded above, then $M$ has a
    semi-injective resolution $M \qra I'$ with $I'_i=0$ for all $i >
    \sup{\H{M}}$.
  \end{prt}
  For an $R$-module $M$, a semi-projective (-injective) resolution is
  just a projective (injective) resolution in the classic sense;
  see~\cite{careil}.

  We use the standard notations $\Ltp{-}{-}$ and $\RHom{-}{-}$ for the
  derived tensor product and derived Hom of complexes; they are
  computed by way of the resolutions described above. Extending the
  usual definitions of Tor and Ext for modules, set
  \begin{equation*}
    \Tor{i}{M}{N} = \H[i]{\Ltp{M}{N}} \qand \Ext{i}{M}{N} =
    \H[-i]{\RHom{M}{N}}
  \end{equation*}
  for $R$-complexes $M$ and $N$ and $i\in\ZZ$. In another extension of
  classic notions, define the projective and injective dimension of an
  $R$-complex by
  \begin{equation*}
    \pd{M} = \inf\setof{\sup{P}}{P \qra M \text{ is a semi-projective
        resolution}}
  \end{equation*}
  and
  \begin{equation*}
    \id{M} = \inf\setof{-\inf{I}}{M \qra I \text{ is a semi-injective
        resolution}}.
  \end{equation*}
\end{bfhpg}

\begin{bfhpg*}[Setup]
  From this point, $R$ denotes a local ring with maximal ideal $\m$
  and residue~field $k = R/\m$. The \emph{embedding dimension} of $R$,
  written $\edim{R}$, is the minimal number of generators of $\m$. The
  \emph{codepth} and \emph{codimension} of $R$ are the differences
  \begin{equation*}
    \codepth{R} = \edim{R} - \dptR \qqand
    \codim{R} = \edim{R} - \dimR,
  \end{equation*}
  where $\dimR$ denotes the Krull dimension of $R$.
\end{bfhpg*}
The depth of an $R$-complex is defined by extension of the homological
characterization of depth of finitely generated modules.

\begin{bfhpg}[Depth]
  \label{finite depth}
  Let $M$ be an $R$-complex. The \emph{depth} of $M$ is defined as
  \begin{equation*}
    \dpt{M} = -\sup{\H{\RHom{k}{M}}}.
  \end{equation*}
  If $\H{M}$ is bounded above, then $M$ has a semi-injective
  resolution $M \qra I$ with $I_i=0$ for $i > \sup{\H{M}}$;
  see \partpgref{res}{I}. Thus, for every $R$-complex $M$ one has
  \begin{equation}
    \label{eq:dptsup}
    \dpt{M} \ge -\sup{\H{M}}.
  \end{equation}
\end{bfhpg}

\begin{bfhpg}[The derived depth formula]
  Let $M$ and $N$ be $R$-complexes. We say that \emph{the derived
    depth formula holds for $M$ and $N$} if there is an equality
  \begin{equation}
    \label{eq:ddf}
    \dpt{\Ltpp{M}{N}} = \dpt{M} + \dpt{N} - \dptR.
  \end{equation}
  Note that this is just a rewrite of the equality \eqref{pdf} in the
  introduction. By \lemcite[(2.1)]{HBF80} the derived depth formula
  holds for complexes $M$ and $N$ if $M$ has finite projective
  dimension and $\H{N}$ is bounded above.
\end{bfhpg}

The next result is due to Dwyer and Greenlees \cite[6.5]{WGDJPG02} and
to Foxby and Iyengar \cite[2.3 and 4.1]{HBFSIn03}.

\begin{prp}
  \protect\pushQED{\qed}%
  \label{prp:supp}
  Let $K$ be the Koszul complex on a set of generators for $\m$. For
  an $R$-complex $M$, the following conditions are equivalent.
  \begin{eqc}
  \item $\H{\Ltp{k}{M}} = 0$;
  \item $\H{\tp{K}{M}} = 0$;
  \item $\H{\Hom{K}{M}} = 0$;
  \item $\H{\RHom{k}{M}} = 0$.\qedhere
  \end{eqc}
\end{prp}

\section{Depth and vanishing of Tate homology---the main theorem}
\label{sec:depth}

\noindent
We start by recalling some facts from \cite{LWCDAJa} and \cite{AIc07}.

\begin{bfhpg}[Complete resolutions]
  An acyclic complex $T$ of projective $R$-modules is called
  \emph{totally acyclic,} if the complex $\Hom{T}{Q}$ is acyclic for
  every projective $R$\nobreakdash-module $Q$. An $R$-module $G$ is
  called \emph{Gorenstein projective} if there exists such a totally
  acyclic complex $T$ with $\Co[0]{T} \is G$.

  Let $M$ be an $R$-complex. A \emph{complete (projective) resolution}
  of $M$ is a diagram
  \begin{equation}
    \label{eq:cpltres}
    T \xra{\tau} P \xra{\pi} M,
  \end{equation}
  where $\pi$ is a semi-projective resolution, $T$ is a totally
  acyclic complex of projective $R$-modules, and $\tau_i$ is an
  isomorphism for $i \gg 0$. The \emph{Gorenstein projective
    dimension} of $M$, written $\Gpd{M}$, is the least integer $n$
  such that there exists a complete resolution \eqref{cpltres} where
  $\tau_i$ is an isomorphism for all $i \ge n$. In particular,
  $\Gpd{M}$ is finite if and only if $M$ has a complete
  resolution. Notice that the homology $\H{M}$ is bounded above if
  $\Gpd{M}$ is finite. Note also that a complex of finite projective
  dimension has finite Gorenstein projective dimension; indeed, $0 \to
  P \to M$ is a complete resolution for every semi-projective
  resolution $P \to M$ with $P$ bounded~above.
\end{bfhpg}

\begin{bfhpg}[Tate homology]
  \label{tate}
  Let $M$ be an $R$-complex of finite Gorenstein projective dimension,
  and let $T \to P \to M$ be a complete resolution. For an $R$-complex
  $N$, the Tate homology of $M$ with coefficients in $N$ is defined as
  \begin{equation*}
    \Ttor{i}{M}{N} =  \H[i]{\tp{T}{N}}.
  \end{equation*}
  This definition is independent of the choice of complete resolution;
  see \cite{AIc07} or \seccite[2]{LWCDAJa} for details. In particular,
  one has
  \begin{equation}
    \label{eq:Ttor=tor}
    \Ttor{i}{M}{N} \is \Tor{i}{M}{N}\quad\text{for } i > \Gpd{M} + \sup{N}.
  \end{equation}
  If $M$ has finite projective dimension or if $N$ is bounded above
  and of finite projective dimension, then $\Ttor{i}{M}{N}=0$ for all
  $i\in\ZZ$; see \prpcite[(2.5) and lem.~(2.7)]{LWCDAJa}.
\end{bfhpg}

The next theorem is our central result. For an $R$-complex $M$ of
finite projective dimension one has $\Ttor{*}{M}{-}=0$, so the theorem
subsumes \lemcite[(2.1)]{HBF80}. The boundedness condition on the
complex $N$, as opposed to its homology, reflects the fact that Tate
homology is not a functor on the derived category over $R$; see the
remarks before \prpcite[(2.5)]{LWCDAJa}.

\begin{thm}
  \label{thm:dpt}
  Let $M$ be an $R$-complex of finite Gorenstein projective dimension
  and let $N$ be a bounded above $R$-complex. If one has
  $\Ttor{i}{M}{N}=0$ for all $i \in \ZZ$, then the derived depth
  formula holds for $M$ and $N$. That is, one has
  \begin{equation*}
    \dpt{\Ltpp{M}{N}} = \dpt{M} + \dpt{N} - \dptR.
  \end{equation*}
\end{thm}

Note that the homology complex $\H{\Ltp{M}{N}}$ is bounded above by
\eqref{Ttor=tor}.

\begin{prf*}
  Choose a complete resolution $T \xra{\tau} P \to M$ and let
  $\qisdef{\pi'}{P'}{N}$ be a semi-projective resolution. The
  quasi-isomorphism $\qisdef{\tp{P}{\pi'}}{\tp{P}{P'}}{\tp{P}{N}}$ is
  then a semi-projective resolution, and the K\"unneth formula yields
  \begin{align*}
    \H{\Ltp{\Ltpp{M}{N}}{k}} &\is \H{\tp{\tpp{P}{P'}}{k}}\\
    &\is \H{\tp[k]{\tpp{P}{k}}{\tpp{P'}{k}}}\\
    &\is \tp[k]{\H{\Ltp{M}{k}}}{\H{\Ltp{N}{k}}}.
  \end{align*}
  It now follows from \eqref{dptsup} and \prpref{supp} that
  $\dpt{\Ltpp{M}{N}}$ is finite if and only if $\dpt{M}$ and $\dpt{N}$
  are both finite. In particular, the left- and right-hand sides of
  the equality we aim to prove are simultaneously finite.

  Assume that both $M$ and $N$ have finite depth.  Consider the
  degreewise split exact sequence of $R$-complexes $0 \to P \to
  \Cone{\tau} \to \Shift{T} \to 0$ and apply the functor $\tp{-}{N}$
  to it. By assumption, the complex $\tp{(\Shift{T})}{N} \is
  \Shift{\tpp{T}{N}}$ is acyclic, so there is a quasi-isomorphism
  $\tp{P}{N} \qra \tp{(\Cone{\tau})}{N}$. In high degrees $K =
  \Cone{\tau}$ is isomorphic to the mapping cone of an isomorphism,
  and the mapping cone of an isomorphism is contractible. Therefore
  there exist homomorphisms $\mapdef{\sigma_i}{K_i}{K_{i+1}}$, such
  that one has $1^{K_i} = \sigma_{i-1}\dif[i]{K} +
  \dif[i+1]{K}\sigma_{i}$ for $i \gg 0$. Since $K$ is a complex of
  projective modules, and $\dif[i+1]{K}\sigma_{i} =
  1^{\Im{\dif[i+1]{K}}}$ holds for $i \gg 0$, it follows that the
  modules $\Ker{\dif[i]{K}} = \Im{\dif[i+1]{K}} \is \Co[i+1]{K}$ are
  projective for $i \gg 0$. Fix $n \gg 0$ and consider the
  contractible subcomplex $J = \cdots \to K_{n+2} \to K_{n+1} \to
  \Im{\dif[n+1]{K}} \to 0$. The sequence $0 \to J \to K \to K/J \to 0$
  is split exact, because the quotient complex $L = K/J = 0 \to
  \Co[n]{K} \to K_{n-1} \to \cdots$ consists of projective modules.
  The complex $\tp{J}{N}$ is contractible, so there are
  quasi-isomorphisms,
  \begin{equation}
    \label{eq:a}
    \tp{P}{N} \qra \tp{K}{N} \qra \tp{L}{N}.
  \end{equation}
  Choose a semi-injective resolution $\mapdef[\qra]{\iota}{N}{I}$,
  where $\iota$ is injective and $I$ is bounded above; see
  \partpgref{res}{I}. Consider the exact sequence $0 \to N \xra{\iota}
  I \to C\to 0$ of $R$\nobreakdash-complexes.  The complex $C =
  \Coker{\iota}$ is bounded above and acyclic, and hence so is the
  complex $\tp{L}{C}$; cf.~\lemcite[2.13]{CFH-06}. Thus, there is a
  quasi-isomorphism
  \begin{equation}
    \label{eq:b}
    \tp{L}{N} \qra \tp{L}{I}.
  \end{equation}
  The complex $\tp{L}{I}$ is bounded above and consists of injective
  $R$\nobreakdash-modules, so it follows from \eqref{a}, \eqref{b},
  and \cite[1.4.I]{LLAHBF91} that it is a semi-injective resolution of
  the complex $\tp{P}{N} \qis \Ltp{M}{N}$. The third equality in the
  next computation follows from \lemref{tev}.
  \begin{align*}
    \dpt{\Ltpp{M}{N}} &= -\sup\H{\RHom{k}{\Ltp{M}{N}}}\\
    &= -\sup\H{\Hom{k}{\tp{I}{L}}}\\
    &= -\sup\H{\tp{\Hom{k}{I}}{L}}\\
    &= -\sup\H{\tp[k]{\Hom{k}{I}}{\tp{k}{L}}}\\
    &= -\sup\H{{\Hom{k}{I}}} - \sup\H{\tp{k}{L}}\\
    &= \dpt{N} - \sup\H{\tp{k}{L}}
  \end{align*}
  For $N=R$ this equality reads
  \begin{equation*}
    \dpt{M} = \dptR  - \sup\H{\tp{k}{L}}.
  \end{equation*}
  The desired equality follows by elimination of the quantity
  $\sup\H{\tp{k}{L}}$.
\end{prf*}

\begin{exa}
  Let $M$ be an $R$-module of finite Gorenstein projective dimension.

  (1) If $N$ is an $R$-module of finite injective dimension, then
  $\Ttor{*}{M}{N}=0$ holds by \lemcite[2.3 and prop.~3.9]{CFH-06}, so
  the derived depth formula holds for $M$ and $N$.

  (2) If $N$ is an $R$-module of finite projective dimension, and
  $\mapdef{\iota}{N}{I}$ is an injective preenvelope, then the derived
  depth formula holds for $M$ and $N' = \Coker{\iota}$, as
  $\Ttor{*}{M}{N}=0$ and $\Ttor{*}{M}{I}=0$ force $\Ttor{*}{M}{N'}=0$.

  (3) If $N$ is an $R$-module of finite injective dimension, and
  $\mapdef{\pi}{P}{N}$ is a projective precover, then the derived
  depth formula holds for $M$ and $\Ker{\pi}$.

  Part (1) is known from \thmcite[6.3]{LWCHHl09}, while (2) and (3)
  appear to be new.
\end{exa}

\begin{bfhpg}[Gorenstein rings]
  \label{gor}
  Let $R$ be Gorenstein. Every $R$-complex with bounded above homology
  has finite Gorenstein projective dimension; see
  \thmcite[3.11]{OVl06}. Thus by \thmref{dpt} the derived depth
  formula holds for $R$-complexes $M$ and $N$ with $\H{M}$ and $N$
  bounded above and $\Ttor{\ast}{M}{N}=0$.

  For finitely generated modules, the Gorenstein projective dimension
  coincides with Auslander and Bridger's notion of G-dimension; see
  \prpcite[3.8]{CFH-06}. The following equality is known as the
  Auslander--Bridger Formula; it holds for every finitely generated
  module $M$ of finite G-dimension,
  \begin{equation}
    \label{eq:ABE}
    \Gdim{M} = \dptR - \dpt{M}.
  \end{equation}
  For finitely generated $R$-modules with $\Ttor{\ast}{M}{N}=0$ there
  is hence an equality
  \begin{equation}
    \label{eq:abedpt}
    \Gdim{\Ltpp{M}{N}} = \Gdim{M} + \Gdim{N},
  \end{equation}
  which represents a natural generalization of \corcite[1.3]{MAs61} to
  G-dimension.
\end{bfhpg}

\begin{rmk}
  \label{rmk:DAJLMS}
  In \cite{DAJLMS04} Jorgensen and \c{S}ega give an example of an
  artinian Gorenstein ring $R$ and a finitely generated $R$-module
  $M$, such that for every $s \ge 0$ there exists a finitely generated
  $R$-module $N_s$ with
  \begin{equation*}
    s = \sup{\H{\Ltp{M}{N_s}}} = -\dpt{\Ltpp{M}{N_s}}.
  \end{equation*}
  Thus, over a Gorenstein ring, boundedness of $\H{\Ltp{M}{N}}$---that
  is, vanishing of $\Tor{\gg 0}{M}{N}$---does not \emph{per se}
  guarantee that the derived depth formula holds. This phenomenon
  disappears over so-called AB rings, where vanishing of homology is
  easier to control.
\end{rmk}

\section{AB rings}
\label{sec:ab}

\noindent
Recall from Huneke and Jorgensen~\cite{CHnDAJ03} that a local ring $R$
is called \emph{AB} if it is Gorenstein, and the following holds for
all finitely generated $R$-modules $M$ and $N$,
\begin{equation*}
  \text{$\Ext{i}{M}{N}
    = 0$ for $i \gg 0$ implies $\Ext{i}{M}{N} = 0$ for $i > \dimR$.}
\end{equation*}
At the end of this section we apply our main theorem to provide a
precise bound for the vanishing of cohomology $\Ext{i}{M}{N}$ for
modules over AB rings; it turns out to depend only on $M$.

\begin{bfhpg}[Tate cohomology]
  \label{TC}
  Let $M$ be an $R$-complex of finite Gorenstein projective dimension,
  and let $T \to P \to M$ be a complete resolution. For an $R$-complex
  $N$, the Tate cohomology of $M$ with coefficients in $N$ is defined
  as
  \begin{equation*}
    \Text{i}{M}{N} =  \H[-i]{\Hom{T}{N}}.
  \end{equation*}
  This definition is independent of the choice of complete resolution;
  see \cite{OVl06} for details. In particular, one has
  \begin{equation}
    \label{eq:Text=ext}
    \Text{i}{M}{N} \is \Ext{i}{M}{N}\quad\text{for } i > \Gpd{M} - \inf{N}.
  \end{equation}
  If $M$ has finite projective dimension or if $N$ is bounded below
  and of finite injective dimension, then one has $\Text{i}{M}{N}=0$
  for all $i\in\ZZ$; see \thmcite[4.5]{OVl06} and
  \lemcite[(4.2)]{LWCDAJa}.
\end{bfhpg}

In some sense, the conditions in the main theorem \thmref[]{dpt} are
easier to verify for moduless over AB rings. Not only is finiteness of
Gorenstein projective dimenion automatic, per the next lemma one only
needs vanishing of homology in high degrees. For modules over a
familiar class of AB rings, namely complete intersections, it becomes
truly easier, as one only needs vanishing of a finite number of
homology modules; see \thmref{ci}.

\begin{prp}
  \label{prp:Ttor-ab}
  Let $R$ be AB and let $M$ and $N$ be finitely generated
  $R$-modules. The following assertions hold.
  \begin{prt}
  \item $\Tor{i}{M}{N}=0$ for $i \gg 0$ implies $\Ttor{i}{M}{N} = 0$
    for all~$i \in \ZZ$.
  \item $\Ext{i}{M}{N} =0$ for $i \gg 0$ implies $\Text{i}{M}{N} =0$
    for all $i\in\ZZ$.
  \end{prt}
\end{prp}

\begin{prf*}
  As $R$ is Gorenstein, $M$ has finite G-dimension. Let $T \to P \to
  M$ be a complete resolution of $M$. For all integers $i$ and $n$
  with $i > n$ one has
  \begin{equation}
    \label{eq:a2}
    \Ttor{i}{M}{N} = \H[i]{\tp{T}{N}} = \Tor{i-n}{\Co[n]{T}}{N}.
  \end{equation}
  The Krull dimension $d = \dimR$ is an upper bound for the
  G-dimension of a finitely generated $R$-module, cf.~\eqref{ABE}. It
  follows from \eqref{Ttor=tor} that there are isomorphisms
  $\Ttor{i}{M}{N} \is \Tor{i}{M}{N}$ for $i > d$. Thus, if the
  homology modules $\Tor{i}{M}{N}$ vanish for $i \gg 0$, then so do
  the modules $\Ttor{i}{M}{N}$. For every $n \in \ZZ$ it follows that
  the modules $\Tor{j}{\Co[n]{T}}{N}$ vanish for $j \gg 0$, and since
  $R$ is AB they vanish for $j>0$; see \thmcite[3.4]{CHnDAJ03}. Now it
  follows from \eqref{a2} that all the Tate homology modules
  $\Ttor{i}{M}{N}$ vanish. The proves part (a); the proof of (b) is
  similar.
\end{prf*}

\begin{rmk}
  \label{rmk:ab}
  If $R$ is AB and $M$ and $N$ are finitely generated $R$-modules with
  $\Tor{i}{M}{N}=0$ for $i \gg 0$, then it follows
  from  \pgref{gor} and \partprpref{Ttor-ab}{a} that the derived depth
  formula holds for $M$ and $N$, and hence that \eqref{abedpt} holds.
\end{rmk}

\begin{rmk}
  Let $R$ be AB and let $M$ and $N$ be finitely generated Gorenstein
  projective $R$-modules with (minimal) complete resolutions $T$ and
  $A$. If one has $\Tor{i}{M}{N}=0$ for $i \gg 0$, then it follows
  from \partprpref{Ttor-ab}{b} and \corcite[(6.2)]{LWCDAJa} that the
  module $\tp{M}{N}$ is Gorenstein projective with (minimal) complete
  resolution $\btp{T}{A}$. Here $\btp{T}{A}$ denotes the pinched
  tensor product defined in \cite{LWCDAJa}.
\end{rmk}

Recall that $R$ is \emph{complete intersection} if there exists a
surjective ring homomorphism $\pi\colon Q \onto \Rhat$, where $Q$ is a
complete regular local ring, and $\Ker{\pi}$ is generated by a
$Q$-regular sequence. The least length of such a sequence equals
$\codim{R}$.  Complete intersection rings are perhaps the best known
examples of AB rings, and they are studied further in the next
section. Here we mention, in passing, a result of Celikbas and Dao
\corcite[1.3]{OClHLD} that constitutes a partial converse to
\rmkref{ab}.

\begin{ipg}
  \label{od}
  Let $R$ be complete intersection of codimension $c$, and assume that
  $R_\p$ is regular for every prime ideal $\p$ in $R$ with $\hgt{\p}
  \le c$.  If the tensor product $\tp{M}{N}$ of two finitely generated
  Gorenstein projective $R$-modules is Gorenstein projective, then one
  has $\Tor{i}{M}{N}=0$ for all $i\ge 1$.
\end{ipg}

The next result on vanishing of cohomology was established for
complete intersection rings by Araya and Yoshino
\thmcite[4.2]{TArYYs98}.

\begin{thm}
  \label{thm:abbound}
  Assume that $R$ is AB and let $M$ and $N$ be finitely generated
  $R$-modules. If one has $\Ext{i}{M}{N}=0$ for $i \gg 0$, then the
  next equality holds,
  \begin{equation*}
    \sup\setof{i \in \ZZ}{\Ext{i}{M}{N}\ne 0} = \dptR - \dpt{M}.
  \end{equation*}
\end{thm}

\begin{prf*}
  As $R$ is Gorenstein, the module $M$ has finite G-dimension. Choose
  a complete resolution $T \xra{\tau} P \to
  M$. By \partprpref{Ttor-ab}{b} the Tate cohomology modules
  $\Text{i}{M}{N}$ vanish for all $i$, so the complex $\Hom{T}{N}$ is
  acyclic.  Let $E$ be the injective hull of the residue field. The
  complex
  \begin{equation*}
    \Hom{\Hom{T}{N}}{E} \is \tp{T}{\Hom{N}{E}},
  \end{equation*}
  is then acyclic as well; the isomorphism is homomorphism evaluation
  \cite[4.4.I]{LLAHBF91} in each degree.  It follows that the Tate
  homology modules $\Ttor{i}{M}{\Hom{N}{E}}$ vanish for all $i$, so
  \thmref{dpt} applies to the modules $M$ and $\Hom{N}{E}$. The latter
  module has depth $0$, so one has
  \begin{equation}
    \label{eq:-dpt}
    \dpt{\Ltpp{M}{\Hom{N}{E}}} = \dpt{M} - \dptR.
  \end{equation}
  The module $\Hom{N}{E}$ is only supported on the maximal ideal of
  $R$, and so are the homology modules of the complex
  $\Ltp{M}{\Hom{N}{E}} \qis \tp{P}{\Hom{N}{E}}$. This explains the
  first equality in the next chain; the second equality is
  homomorphism evaluation.
  \begin{equation}
    \label{eq:compu}
    \begin{split}
      -\dpt{\tpp{P}{\Hom{N}{E}}} &= \sup\H{\tp{P}{\Hom{N}{E}}}\\
      &= \sup\H{\Hom{\Hom{P}{N}}{E}}\\
      &= -\inf{\H{\Hom{P}{N}}}\\%
      &= -\inf{\H{\RHom{M}{N}}}\\
      &= \sup\setof{i \in \ZZ}{\Ext{i}{M}{N}\ne 0}
    \end{split}
  \end{equation}
  The desired equality now follows from \eqref{-dpt} and
  \eqref{compu}.
\end{prf*}

\section{Complete intersections}
\label{sec:ci}

\noindent
Homology of finitely generated modules over complete intersection
rings is rigid in the following sense; see \cite[remarks before
thm.~1.9]{HJW-01}.

\begin{bfhpg}[Fact]
  \label{rigid}
  Let $R$ be complete intersection of codimension $c$ and let $M$ and
  $N$ be finitely generated $R$-modules. Let $n \ge 0$ be an integer;
  if one has $\Tor{i}{M}{N}=0$ for all $i$ with $n+c \ge i \ge n$,
  then one has $\Tor{i}{M}{N}=0$ for all $i \ge n$.
\end{bfhpg}

Combined with \thmref{dpt} and \thmcite[4.9]{LLAROB00} this fact has the
following consequence.

\begin{cor}
  \label{cor:rigid}
  Let $R$ be complete intersection of codimension $c$, and let $M$ and
  $N$ be finitely generated $R$-modules. If one has $\Tor{i}{M}{N}=0$
  for $c+1$ consecutive values of $i \ge 0$, then the Tate homology
  modules $\Ttor{i}{M}{N}$ vanish for all $i\in\ZZ$, and the derived
  depth formula holds for $M$ and $N$. \qed
\end{cor}

One goal of this section is to obtain a similar result for modules
that are not finitely generated; see \thmref{ci}.

Recall from \cite{AGP-97} that a (\emph{codimension} $c$)
\emph{quasi-deformation} of $R$ is a diagram of local homomorphism $R
\xra{\rho} R' \xla{\pi} Q$, where $\rho$ is flat, and $\pi$ is
surjective with kernel generated by a $Q$-regular sequence (of length
$c$).

\begin{lem}
  \label{lem:ci}
  Let $M$ and $N$ be $R$-complexes with $N$ bounded above.  If there
  exists a codimension $c$ quasi-deformation $R \to R' \leftarrow Q$
  such that $\pd[Q]{\tpp{R'}{M}}$ is finite, then $\Gpd{M}$ is finite,
  and the following conditions are equivalent.
  \begin{eqc}
  \item $\Ttor[R]{i}{M}{N}=0$ for all $i\in\mathbb Z$.
  \item $\Tor{i}{M}{N}=0$ for all $i\gg 0$.
  \item[\eqclbl{ii'}] $\Ttor{i}{M}{N}=0$ for all $i \ll 0$.
  \item $\Tor{i}{M}{N}=0$ for $c+1$ consecutive values of $i> \Gpd M +
    \sup{N}$.
  \item $\Ttor[R]{i}{M}{N}=0$ for $c+1$ consecutive values of $i$.
  \end{eqc}
\end{lem}

\noindent To not interrupt the flow, we defer the proof of this lemma
to the end of the section.

We say that an $R$-complex $M$ has \emph{finite CI-dimension}, if
there is a quasi-deformation of $R$ such that $\pd[Q]{\tpp{R'}{M}}$ is
finite.  For modules $M$ and $N$ with $\Tor{\ge 1}{M}{N} =0$ the next
theorem recovers Iyengar's \thmcite[4.3]{SIn99}; it also subsumes a
recent generalization of this result due to Sahandi, Sharif, and
Yassemi \thmcite[3.3]{SSY-}.

\begin{thm}
  \label{thm:cim}
  Let $M$ and $N$ be $R$-complexes with $\H{N}$ bounded above. If $M$
  has finite CI-dimension, and one of the following conditions holds.
  \begin{prt}
  \item one has $\Tor{i}{M}{N}=0$ for $i \gg 0$; or
  \item the complex $N$ is bounded above, and one has
    $\Ttor{i}{M}{N}=0$ for $i \ll 0$;
  \end{prt}
  then the Tate homology modules $\Ttor{i}{M}{N}$ vanish for all
  $i\in\ZZ$, and the derived depth formula holds for $M$ and $N$.
\end{thm}

We precede the proof with a technical observation.

\begin{rmk}
  \label{rmk:trunc}
  Let $N$ be an $R$-complex with $\H{N}$ bounded above. Set $s =
  \sup{\H{N}}$ and let $N'$ be the soft truncation of $N$ at $s$,
  i.e.\ the bounded above complex $0 \to \Co[s]{N} \to N_{s-1} \to
  \cdots$. The natural morphism \mbox{$N \onto N'$} is a
  quasi-isomorphism; in particular, one has $\dpt{N'} = \dpt{N}$. For
  every $R$-complex $M$ there are isomorphisms
  \begin{equation*}
    \Tor{i}{M}{N'} \is \Tor{i}{M}{N} \ \text{ for all $i\in\ZZ$}.
  \end{equation*}
  Moreover, one has $\dpt{\Ltpp{M}{N'}} =\dpt{\Ltpp{M}{N}}$, so the
  depth formula holds for $M$ and $N$ if and only if it holds for $M$
  and $N'$.
\end{rmk}

\begin{bfhpg*}[Proof of Theorem \pgref{thm:cim}]\protect\pushQED{\qed}
  Under the assumption that $N$ is bounded above and
  $\Ttor{i}{M}{N}=0$ holds for $i \ll 0$, the assertions follow immediately from
  \lemref{ci} and \thmref{dpt}.

  Assume now that $\Tor{i}{M}{N}=0$ holds for $i \gg 0$. By
  \rmkref{trunc} we can replace $N$ by its soft truncation at
  $\sup{\H{N}}$; that is, we can assume that $N$ is bounded above. Now
  the assertions follow as above.\qed
\end{bfhpg*}

Over a complete intersection ring, vanishing of Tate homology
$\Ttor{i}{M}{N}$ for all $i \in \ZZ$ can be inferred from a finite gap
in (Tate) homology, and the length of that gap is independent of $M$
and $N$.

\begin{thm}
  \label{thm:ci}
  Let $R$ be complete intersection of codimension $c$ and
  dimension~$d$. Let $M$ and $N$ be $R$-complexes with $\H{M}$ and
  $\H{N}$ bounded above. If one of the following conditions holds,
  \begin{prt}
  \item $\Tor{i}{M}{N}=0$ for $c+1$ consecutive values of $i > d +
    \sup{\H{M}} + \sup{\H{N}}$;~or
  \item the complex $N$ is bounded above, and one has
    $\Ttor{i}{M}{N}=0$ for $c+1$ consecutive values of $i$;
  \end{prt}
  then the Tate homology modules $\Ttor{i}{M}{N}$ vanish for all
  $i\in\ZZ$, and the derived depth formula holds for $M$ and $N$.
\end{thm}

\begin{prf*}
  A complete intersection ring is Gorenstein, so by
  \thmcite[3.10]{OVl06} one has
  \begin{equation*}
    \Gpd{M} \le d + \sup{\H{M}} < \infty.
  \end{equation*}
  By assumption there is a homomorphism $\pi\colon Q \onto \Rhat$,
  where $Q$ is a complete regular local ring, and $\Ker{\pi}$ is
  generated by a $Q$-regular sequence $x_1,\ldots,x_c$. Thus, the
  diagram $R \to \Rhat \twoheadleftarrow Q$ is a codimension $c$
  quasi-deformation of $R$, and because $\H{\tp{\Rhat}{M}}$ is bounded
  above, one has $\pd[Q]{\tpp{\Rhat}{M}} < \infty$.

  If $N$ is bounded above and one has $\Ttor{i}{M}{N}=0$ for $c+1$
  consecutive values of $i$, then \lemref{ci} yields
  $\Ttor{i}{M}{N}=0$ for all $i\in\ZZ$.  As $\Gpd{M}$ is finite, the
  derived depth formula holds for $M$ and $N$ by \thmref{dpt}.

  Now, set $s = \sup{\H{N}}$ and assume that one has $\Tor{i}{M}{N}=0$
  for $c+1$ consecutive values of $i > d+ \sup{\H{M}} + s$. By
  \rmkref{trunc} we can replace $N$ by its soft truncation at $s$;
  that is, we can assume that $N$ is bounded above with $\sup{N} =
  s$. Now the assertions follow as above.
\end{prf*}

\begin{rmk}
  Over a Gorenstein ring $R$ of dimension $d$, the number
  $d+\sup{\H{M}}$ is an upper bound for the Gorenstein projective
  dimension of every $R$-complex $M$. The proof above and the fact
  that Tate homology is balanced, see \seccite[5]{LWCDAJa}, shows that
  one can replace the quantity $d + \sup{\H{M}} + \sup{\H{N}}$ in the
  theorem with $\min\set{\Gpd{M}+\sup{\H{N}}, \Gpd{N}+\sup{\H{M}}}$.
\end{rmk}

For complexes with bounded and degreewise finitely generated homology,
for finitely generated modules in particular, the CI-dimension agrees
with Avramov, Gasharov, and Peeva's notion of CI-dimension; see
\cite{AGP-97, SSW04}.

\begin{cor}
  \label{cor:cidim}
  Let $R$ be complete intersection of codimension $c$ and let $M$ and
  $N$ be finitely generated $R$-modules. If one has $\Ttor{i}{M}{N}=0$
  for $c+1$ consecutive values of $i$ or $\Tor{i}{M}{N}=0$ for $c+1$
  consecutive values of $i \ge 0$, then the following equality holds,
  \begin{equation*}
    \CIdim{\Ltpp{M}{N}} = \CIdim{M} + \CIdim{N}.
  \end{equation*}
\end{cor}

\begin{prf*}
  It follows from the theorem and \corref{rigid} that the derived
  depth formula holds for $M$ and $N$.  For every $R$-complex $X$ with
  $\H{X}$ finitely generated, there is an Auslander--Buchsbaum-type
  formula $\CIdim{X} = \dptR - \dpt{X}$; see \prpcite[3.3]{SSW04}. Now
  the desired equality follows from the depth formula.
\end{prf*}

\begin{thm}
  \label{thm:cif}
  Let $R$ be of codepth $c$, let $M$ be a finitely generated
  $R$-module of finite CI-dimension, and let $N$ be an $R$-module. If
  one has $\Ttor{i}{M}{N}=0$ for $c+1$ consecutive values of $i$ or
  $\Tor{i}{M}{N}=0$ for $c+1$ consecutive values of $i > \dptR$, then
  the Tate homology modules $\Ttor{i}{M}{N}$ vanish for all $i\in\ZZ$,
  and the derived depth formula holds for $M$~and~$N$.
\end{thm}

\begin{prf*}
  By \thmcite[(1.4) and (5.6)]{AGP-97} the number $d = \dptR$ is an
  upper bound for the CI-dimension of $M$, and the complexity of $M$
  is at most $c$. If one has $\Tor{i}{M}{N}=0$ for $c+1$ consecutive
  values of $i > d$, then \corcite[2.3]{DAJ01} yields
  $\Tor{i}{M}{N}=0$ for all $i > d$, and then it follows from
  \lemref{ci} and \thmref{dpt} that the derived depth formula holds
  for $M$ and $N$.

  Let $T \to P \to M$ be a complete resolution. By
  \lemcite[(1.5)]{AGP-97} every syzygy of $M$ has finite
  CI-dimension. In particular, $\Co[d]{T}$ has finite CI-dimension,
  and it follows that $\Co[i]{T}$ has finite CI-dimension for every
  $i\in\ZZ$. If one has $\Ttor{i}{M}{N}=0$ for $c+1$ consecutive
  values of $i$, then there exists an $\imath \in \ZZ$ such that one
  has
  \begin{equation*}
    0 =
    \Tor{d+1}{\Co[\imath]{T}}{N}=\Tor{d+2}{\Co[\imath]{T}}{N}=\cdots =
    \Tor{d+c+1}{\Co[\imath]{T}}{N}.
  \end{equation*}
  Now \corcite[2.3]{DAJ01} yields $\Tor{i}{\Co[\imath]{T}}{N}=0$ for
  all $i > d$, and then it follows from \lemref{ci} and \thmref{dpt}
  that the desired formula holds for $M$ and $N$.
\end{prf*}

\begin{bfhpg*}[Proof of Lemma \pgref{lem:ci}]\protect\pushQED{\qed}    
  \setcounter{equation}{0}%
  \renewcommand{\theequation}{\arabic{equation}} The complex $M$ has
  finite CI-dimension. If the homology complex $\H{M}$ is bounded,
  then it follows from \thmcite[5.1.(b) and rmk.~2.5]{SSW08} that
  $\Gpd{M}$ is finite. Here we give a direct but similar argument that
  does not use boundedness of $\H{M}$.

  Let $R \xra{\rho} R' \xla{\pi} Q$ be a quasi-deformation such that
  $\pd[Q]{\tpp{R'}{M}}$ is finite. Assume, without loss of generality,
  that the rings $Q$ and, therefore, $R'$ are complete. Then $Q$ has a
  dualizing complex $D$, and the complex $D' = \RHom[Q]{R'}{D}$ is
  dualizing for $R'$; see \seccite[3]{SInHKr06}.  The complex
  $\H{\Ltp[Q]{D}{\tpp{R'}{M}}}$ is bounded above by \thmcite[2.4.F and
  2.3.F]{LLAHBF91}, and it follows from \prpcite[7.3]{SInHKr06} that
  the complex $M$ belongs to the Auslander category
  $\hat{\mathcal{A}}(R)$. Let $x_1,\ldots,x_c$ be a $Q$-regular
  sequence that generates $\Ker{\pi}$. The Koszul complex on
  $x_1,\ldots,x_c$ is a projective resolution of $R'$ over $Q$, so one
  has $D' \qis \Shift[-c]{\Ltpp[Q]{D}{R'}}$. Now it is straightforward
  to verify that $\H{\Ltp[R']{D'}{\tpp{R'}{M}}}$ is bounded above and
  that $\tp{R'}{M}$ belongs to $\hat{\mathcal{A}}(R')$. Therefore,
  $g=\Gpd[R']{\tpp{R'}{M}}$ is finite by \thmcite[8.1]{SInHKr06}. Let
  \mbox{$P \qra M$} be a semi-projective resolution over $R$. To prove
  that $\Gpd{M}$ is finite, it suffices to show that the module
  $\Co[g]{P}$ has finite Gorenstein projective dimension or,
  equivalently, finite Gorenstein flat dimension;
  cf.~\thmcite[4.1]{CFH-06}. The quasi-isomorphism $\tp{R'}{P} \to
  \tp{R'}{M}$ is a semi-projective resolution over $R'$. By
  \thmcite[3.4]{OVl06} the module $\Co[g]{\tp{R'}{P}} \is
  \tp{R'}{\Co[g]{P}}$ is Gorenstein projective, in particular it is
  Gorenstein flat, and then so is $\Co[g]{P}$ by
  \thmcite[A]{LWCSSW10}.  Thus $\Gpd{M}$ is finite.

  The implications
  \begin{equation*}
    (iv) \impliedby (ii') \impliedby (i)\implies(ii)\implies(iii)\implies(iv)
  \end{equation*}
  are trivial or follow from \eqref{Ttor=tor}.  It remains to prove
  that $(iv)$ implies $(i)$.

  Let $T\to P\to M$ be a complete resolution. The morphism $\tp{R'}{P}
  \to \tp{R'}{M}$ is then a semi-projective resolution over $R'$, and
  it is straightforward to verify that $\tp{R'}{T}$ is a totally
  acyclic complex of projective $R'$-modules. Thus, by flatness of
  $R'$ one has
  \begin{equation*}
    \tp{R'}{\Ttor{i}{M}{N}}
    \is \Ttor[R']{i}{\tp{R'}{M}}{\tp{R'}{N}}
  \end{equation*}
  for all $\iinZ$. Without loss of generality, assume that $\rho$ is
  the identity map. The assumptions are now that $R$ is isomorphic to
  $Q/(x_1,\ldots,x_c)$, that $\pd[Q]{M}$ is finite, and that there is
  an integer $h$ such that
  \begin{equation}
    \label{eq:ttors}
    0 = \Ttor[R]{h}{M}{N}=\Ttor[R]{h+1}{M}{N}= \cdots
    =\Ttor[R]{h+c}{M}{N}.
  \end{equation}
  For $j\in\ZZ$ set $M_j=\Co[j]{T}$. Notice that because $R$ and $M$
  have finite projective dimension over $Q$, each module $M_j$ has
  finite projective dimension over $Q$ as well. Set $d = \dpt[]{Q}$;
  fix a $j$ and let $F$ be a projective resolution of $M_j$ over $Q$
  of length $\pd[Q]{M_j} \le d$. Let $L$ be a semi-projective
  resolution of $N$ over $R$.  Without loss of generality, assume that
  one has $\sup{N}=0$. The filtrations $\mathcal{F}$ and $\mathcal{L}$
  defined by
  \begin{equation*}
    (\mathcal{F}^p \tpp[Q]{F}{L})_n = \coprod_{i\le p}F_i\otimes_Q L_{n-i}
    \qand
    (\mathcal{L}^p\tpp[Q]{L}{F})_n=\coprod_{i\le p}L_{i}\otimes_Q F_{n-i}
  \end{equation*}
  are bounded; they give rise to spectral sequences
  \begin{equation*}
    \Eterm{\mathcal{F}}{p}{q} \underset p \implies
    \H[n]{\tp[Q]{F}{L}} \qqand \Eterm{\mathcal{L}}{p}{q}
    \underset p \implies \H[n]{\tp[Q]{L}{F}}.
  \end{equation*}
  The $\mathrm{E}^2$-terms are the iterated homologies of the
  underlying double complexes, obtained by first taking homology along
  columns and then along rows. The terms
  \begin{equation*}
    \Eterm{\mathcal{F}}{p}{q} = \Tor[Q]{p}{M_j}{\H[q]{N}}
  \end{equation*}
  vanish for $q > 0$, and they vanish for $p$ not in
  $\set{0,\ldots,d}$. In particular, for $n > d$ one has
  $\Eterm{\mathcal{F}}{i}{n-i}=0$ for all $i\in\ZZ$. It follows that
  the homology modules $\H[n]{\tp[Q]{F}{L}}$ vanish for $n > d$; see
  \prpcite[5.5]{careil}. From the isomorphism $\H{\tp{F}{L}} \is
  \H{\tp{L}{F}}$ and \prpcite[5.3a]{careil} one now gets
  \begin{equation}
    \label{eq:infterm}
    \Einfterm{\mathcal{L}}{i}{n-i}=0 \ \text{ for all $n>d$ and all $i\in\ZZ$.}
  \end{equation}

  The homology in degree $q$ within the $p$th column in the double
  complex $\tp{L}{F}$ is isomorphic to
  $\tp{L_p}{\H[q]{\Ltp[Q]{R}{M_j}}}$.  The Koszul complex
  $\Koszul[Q]{x_1,\ldots,x_c}$ is a free resolution of $R$ over $Q$,
  and the elements $x_1,\ldots,x_c$ act trivially on $M_j$, so one has
  $\tp{L_p}{\H[q]{\Ltp[Q]{M_j}{R}}} \is
  \tp{L_p}{{M_j}^{\binom{c}{q}}}$ for all $q$.  Thus the
  $\mathrm{E}^2$-terms in the second sequence are
  \begin{equation*}
    \Eterm{\mathcal{L}}{p}{q} = \Tor{p}{N}{M_j^{\binom{c}{q}}} \is
    \Tor{p}{M_j}{N}^{\binom{c}{q}}.
  \end{equation*}
  Clearly the terms $\Eterm{\mathcal{L}}{p}{q}$ vanish for $q$ not in
  $\set{0,\ldots,c}$, and $\Eterm{\mathcal{L}}{p}{q}$ vanishes for all
  $q$ in $\set{0,\ldots,c}$ if and only if it vanishes for one of
  them.

  Assume now that $j < h-d+c-1$ holds. One then has $h-j > 0$, and the
  hypothesis \eqref{ttors} yields $\Eterm{\mathcal{L}}{p}{q} = 0$ for
  $p$ in $\set{h-j, \ldots, h-j +c}$. This yields the limit terms
  \begin{equation}
    \label{eq:limterm}
    \Eterm{\mathcal{L}}{h-j+c+1}{0} = \Einfterm{\mathcal{L}}{h-j+c+1}{0} \qqand
    \Eterm{\mathcal{L}}{h-j-1}{c} = \Einfterm{\mathcal{L}}{h-j-1}{c}.
  \end{equation}
  Moreover, the inequalities $h-j+c+1 > d$ and $c + (h-j-1) >d$ hold,
  so \eqref{infterm} and \eqref{limterm} combine to yield
  $\Eterm{\mathcal{L}}{h-j+c+1}{0}=0$ and
  $\Eterm{\mathcal{L}}{h-j-1}{c} =0$ and, therefore,
  $\Eterm{\mathcal{L}}{h-j+c+1}{q}=0=\Eterm{\mathcal{L}}{h-j-1}{q}$
  for all $q$.

  Iterating this argument, one gets $\Eterm{\mathcal{L}}{p}{q}=0$ for
  all $p>0$ and all $q$. Thus, for every $j < h-d+c-1$ one has $0 =
  \Tor{p}{M_j}{N} = \Ttor{p+j}{M}{N}$ for all $p>0$; hence
  $\Ttor{i}{M}{N}=0$ holds for all $i \in \ZZ$.\qed
\end{bfhpg*}

\section{Auslander's depth formula}
\label{sec:formula}

\noindent
The purpose of this section is to connect the derived depth formula
\eqref{ddf} with another generalization, due to Auslander
\cite{MAs61}, of the Auslander--Buchsbaum Formula.

In view of \lemref{ci} the next result generalizes
\thmcite[3]{SChSIn01} by replacing finitely generated modules with
complexes of modules with mild boundedness conditions and no
assumptions of finite generation.

\begin{thm}
  \label{thm:s}
  Let $M$ be an $R$-complex of finite Gorenstein projective dimension
  and let $N$ be a bounded above $R$-complex. If one has
  $\Ttor{i}{M}{N}=0$ for all $i \in \ZZ$, then $s =
  \sup{\H{\Ltp{M}{N}}}$ is finite, there is an inequality,
  \begin{equation*}
    - s \le \dpt{M} + \dpt{N} - \dptR,
  \end{equation*}
  and equality holds if and only if $\dpt{\H[s]{\Ltp{M}{N}}}$ is zero.
  Moreover,~if~one has $s=0$ or $\dpt{\H[s]{\Ltp{M}{N}}} \le 1$, then
  there is an~equality
  \begin{equation*}
    \dpt{\H[s]{\Ltp{M}{N}}} - s = \dpt{M} + \dpt{N} - \dptR.
  \end{equation*}
\end{thm}

\begin{prf*}
  By \thmref{dpt} the complex $\H{\Ltp{M}{N}}$ is bounded above, and
  one has $\dpt{\Ltpp{M}{N}} = \dpt{M} + \dpt{N} - \dptR$; the
  inequality now follows from \eqref{dptsup}. By
  \cite[1.5.(3)]{HBFSIn03} there is an isomorphism
  \begin{equation*}
    \H[-s]{\RHom{k}{\Ltp{M}{N}}} \is \Hom{k}{\H[s]{\Ltp{M}{N}}},
  \end{equation*}
  and it follows that the equality $\dpt{\Ltpp{M}{N}} = -s$ holds if
  an only if the module $\H[s]{\Ltp{M}{N}}$ has depth zero. Finally,
  by \thmcite[2.3]{SIn99} the desired equality $\dpt{\Ltpp{M}{N}} =
  \dpt{\H[s]{\Ltp{M}{N}}} - s$ holds provided that one has
  $\dpt{\H[s]{\Ltp{M}{N}}} - s \le \dpt{\H[i]{\Ltp{M}{N}}} - i$ for
  all $i \ge s$.
\end{prf*}

\begin{bfhpg}[The depth formula]
  Let $M$ and $N$ be finitely generated $R$-modules. Following Choi
  and Iyengar~\cite{SChSIn01} we say that the \emph{depth formula
    holds for $M$ and $N$} if $s = \sup\setof{i}{\Tor{i}{M}{N} \ne 0}$
  is finite and one has
  \begin{equation*}
    \dpt{\Tor{s}{M}{N}}
    - s = \dpt{M} + \dpt{N} - \dptR.
  \end{equation*}
  Auslander \thmcite[1.2]{MAs61} proved that the formula holds, if $M$
  has finite projective dimension, and one has $s=0$ or
  $\dpt{\Tor{s}{M}{N}} \le 1$.
\end{bfhpg}

\begin{cor}
  Let $M$ and $N$ be finitely generated $R$-modules, and assume that $s
  = 0$ or $\dpt{\Tor{s}{M}{N}} \le 1$. The depth formula holds for $M$
  and $N$ if one of the following conditions is satisfied.
  \begin{prt}
  \item $M$ has finite G-dimension, and one has
    $\Ttor{i}{M}{N}=0$ for all $i \in \ZZ$.
  \item $R$ is AB, and one has $\Tor{i}{M}{N}=0$ for $i \gg 0$.
  \item $R$ is complete intersection of codimension $c$, and one has
    $\Ttor{i}{M}{N}=0$ for $c+1$ consecutive values of $i$.
  \item $R$ is complete intersection of codimension $c$, and one has
  $\Tor{i}{M}{N}=0$ for $c+1$ consecutive values of $i \ge 0$.
\item $M$ has finite CI-dimension, and one has
  $\Ttor{i}{M}{N}=0$ for $\codepth{R}+1$ consecutive values of $i\ge 0$.
\item $M$ has finite CI-dimension, and one has
  $\Tor{i}{M}{N}=0$ for $\codepth{R}+1$ consecutive values of $i > \dptR$.
  \end{prt}
\end{cor}

\begin{prf*}
  Part (a) is immediate from \thmref{s}. Part (b) follows from (a) in
  view of \partprpref{Ttor-ab}{a}. Similarly, part (c) follows in view
  of \thmref{ci}, part (d) follows in view of \corref{rigid}, and
  parts (e) and (f) follow in view of \thmref{cif}.
\end{prf*}

\section{Vanishing of cohomology}
\label{sec:coh}

\noindent
Let $M$ and $N$ be finitely generated $R$-modules. If $M$ has finite
projective dimension or $N$ has finite injective dimension, then the
largest index for which $\Ext{i}{M}{N}$ does not vanish is $i = \dptR
- \dpt{M}$. That is, the vanishing of cohomology $\Ext{i}{M}{N}$ only
depends on $M$; see Ischebeck~\cite[2.6]{FIs69}. The~next~result is
more general, as Tate cohomology $\Text{*}{M}{N}$ vanishes if $M$ has
finite projective dimension or $N$ has finite injective dimension, see
\pgref{TC}, and finite projective/injective dimension implies finite
G-dimension/Gorenstein injective dimension.

\begin{thm}
  \label{thm:appl}
  Let $M$ and $N$ be finitely generated $R$-modules such that $M$ has
  finite $G$-dimension or $N$ has finite Gorenstein injective
  dimension. If one has $\Text{i}{M}{N}=0$ for all $i \in \ZZ$, then
  the next equality holds
  \begin{equation*}
    \sup\setof{i \in \ZZ}{\Ext{i}{M}{N}\ne 0} = \dptR - \dpt{M}.
  \end{equation*}
\end{thm}
\noindent
The proof is given at the end of the section; note that under the
assumption that $M$ has finite G-dimension, the desired equality
follows from the proof of \thmref{abbound}.

The notion of Gorenstein injective dimension is dual to that of
Gorenstein projective dimension, and Tate cohomology $\Text{*}{M}{N}$
can be extended to the situation where the second variable $N$ has
finite Gorenstein injective dimension; see \cite{JAsSSl06,LWCDAJa}.

Before we start the proof of \thmref{appl} we record an easy
consequence.

\begin{rmk}
  Let $M$ and $N$ be finitely generated $R$-modules. Under each of the
  following conditions,
  \begin{prt}
  \item $M$ has finite $G$-dimension and $\pd{N}$ is finite
  \item $N$ has finite Gorenstein injective dimension and $\id{M}$ is
    finite
  \end{prt}
  Tate cohomology $\Text{*}{M}{N}$ vanishes by \thmcite[4.5]{OVl06}
  and \thmcite[3.9]{JAsSSl06}, whence one has $\sup\setof{i \in \ZZ}{\Ext{i}{M}{N}\ne 0} = \dptR - \dpt{M}$.
\end{rmk}

As a first step towards a proof of \thmref[]{appl} we recall the
notion of width.

\begin{bfhpg}[Width]
  The \emph{width} of an $R$-complex $M$ is defined as
  \begin{equation*}
    \wdt{M} = \inf\H{\Ltp{k}{M}}.
  \end{equation*}
  There is an obvious inequality
  \begin{equation}
    \label{eq:wdtinf}
    \wdt{M} \ge \inf\H{M},
  \end{equation}
  and equality holds if $\H{M}$ is bounded below and degreewise
  finitely generated.
\end{bfhpg}

\begin{prp}
  \label{prp:wdtI}
  Let $N$ be an $R$-complex of finite Gorenstein injective dimension
  and let $M$ be a bounded above $R$-complex. If one has
  $\Text{i}{M}{N}=0$ for all $i \in \ZZ$, then the next equality
  holds,
  \begin{equation*}
    \wdt{\RHom{M}{N}} = \dpt{M} + \wdt{N} - \dptR.
  \end{equation*}
\end{prp}

By \cite[(5.6.1)]{LWCDAJa} the homology complex $\H{\RHom{M}{N}}$ is
bounded~\mbox{below.}

\begin{prf*}
  Choose a complete injective resolution $N \to I \xra{\upsilon} U$
  and let $\qisdef{\pi}{P}{M}$ be a semi-projective resolution. The
  induced quasi-isomorphism
  \begin{equation*}
    \qisdef{\Hom{\pi}{I}}{\Hom{M}{I}}{\Hom{P}{I}}
  \end{equation*}
  is a semi-injective resolution, and the K\"unneth formula yields
  \begin{align*}
    \H{\RHom{k}{\RHom{M}{N}}} &\is \H{\Hom{k}{\Hom{P}{I}}}\\
    &\is \H{\Hom{\tp{k}{P}}{I}}\\
    &\is \H{\Hom[k]{\tp{k}{P}}{\Hom{k}{I}}}\\
    &\is \Hom[k]{\H{\Ltp{k}{M}}}{\H{\RHom{k}{N}}}.
  \end{align*}
  It follows from \eqref{wdtinf} and \prpref{supp} that
  $\wdt{\RHom{M}{N}}$ is finite if and only if $\wdt{M}$ and $\dpt{N}$
  are both finite. In particular, the left- and right-hand sides of
  the desired equality are simultaneously finite.

  Assume that $\wdt{M}$ and $\dpt{N}$ are finite. Set $K =
  \Shift[-1]{(\Cone{\iota})}$; consider the degreewise split exact
  sequence $0 \to \Shift[-1]{U} \to K \to I \to 0$ and apply the
  functor $\Hom{M}{-}$. By assumption, the complex $\Hom{M}{U}$ is
  acyclic, cf.~\dfncite[(5.5)]{LWCDAJa}, so there is a quasi
  isomorphism $\Hom{M}{K} \qra \Hom{M}{I}$.  In low degrees, $K$ is
  isomorphic to the mapping cone of an isomorphism. Therefore, there
  exist homomorphisms $\mapdef{\sigma_i}{K_i}{K_{i+1}}$ such that
  $1^{K_i} = \sigma_{i-1}\dif[i]{K} + \dif[i+1]{K}\sigma_{i}$ holds
  for $i \ll 0$. Since $K$ is a complex of injective modules, and one
  has $\dif[i+1]{K}\sigma_{i} = 1^{\Im{\dif[i+1]{K}}}$ for $i \ll 0$,
  it follows that the modules $\Ker{\dif[i]{K}} = \Im{\dif[i+1]{K}}
  \is \Co[i+1]{K}$ are injective for $i \ll 0$. Fix $n \ll 0$; the
  subcomplex $J = \cdots \to K_{n+2} \to K_{n+1} \to \Im{\dif[n+1]{K}}
  \to 0$ consists of injective modules, the sequence $0 \to J \to K
  \to K/J \to 0$ is split exact, and the quotient complex $K/J = 0 \to
  \Co[n]{K} \to K_{n-1} \to \cdots$ is contractible.  It follows that
  there are quasi-isomorphisms,
  \begin{equation}
    \label{eq:1a}
    \Hom{M}{J} \qra \Hom{M}{K} \qra \Hom{M}{I}.
  \end{equation}

  Choose a semi-injective resolution $\mapdef[\qra]{\iota'}{M}{I'}$,
  where $\iota'$ is injective and $I'$ is bounded above; see
  \partpgref{res}{I}. Then the complex $C = \Coker{\iota'}$ is bounded
  above and acyclic. Apply $\Hom{-}{J}$ to the exact sequence $0 \to M
  \xra{\smash{\iota'}} I' \to C \to 0$ of $R$-complexes.  The complex
  $\Hom{C}{J}$ is acyclic by \lemcite[2.5]{CFH-06}, so there is a
  quasi-isomorphism
  \begin{equation}
    \label{eq:1b}
    \Hom{I'}{J} \qra \Hom{M}{J}.
  \end{equation}
  The complex $\Hom{I'}{J}$ is bounded below and consists of flat
  $R$-modules, so it is semi-flat in the sense of \cite{LLAHBF91,dga},
  and it follows from \eqref{1a} and \eqref{1b} that there is a
  quasi-isomorphism $\Hom{I'}{J} \qra \Hom{M}{I}$. The third equality
  in the next chain uses homomorphism evaluation, a variation on
  \lemcite[4.4.I]{LLAHBF91} which is proved similarly to \lemref{tev}.
  \begin{align*}
    \wdt{\RHom{M}{N}} &= \inf\H{\Ltp{k}{\RHom{M}{N}}}\\
    &= \inf\H{\tp{k}{\Hom{I'}{J}}}\\
    &= \inf\H{\Hom{\Hom{k}{I'}}{J}}\\
    &= \inf\H{\Hom[k]{\Hom{k}{I'}}{\Hom{k}{J}}}\\
    &= \inf\H{\Hom{k}{J}} - \sup\H{\Hom{k}{I'}}\\
    &= \inf\H{\Hom{k}{J}} + \dpt{M}.
  \end{align*}
  For $M=R$ this equality reads $\wdt{N} = \inf\H{\Hom{k}{J}} +
  \dpt{R}$, and the desired equality follows.
\end{prf*}

A dual argument yields the next result, which is also invoked in the
proof of \thmref{appl}. The special case of \prpref{wdtP} where $M$
and $N$ are finitely generated, which is the context of
\thmref[]{appl}, follows from the proof of \thmref{abbound}.

\begin{prp}
  \protect\pushQED{\qed}%
  \label{prp:wdtP}
  Let $M$ be an $R$-complex of finite Gorenstein projective dimension
  and let $N$ be a bounded below $R$-complex. If one has
  $\Text{i}{M}{N}=0$ for all $i \in \ZZ$, then the next equality
  holds,
  \begin{equation*}
    \wdt{\RHom{M}{N}} = \dpt{M} + \wdt{N} - \dptR.\qedhere
  \end{equation*}
\end{prp}

Note that by \eqref{Text=ext} the homology complex $\H{\RHom{M}{N}}$
is bounded~\mbox{below.}

\begin{bfhpg*}[Proof of Theorem \pgref{thm:appl}]
  \protect\pushQED{\qed}%
  It follows from \cite[(5.6.1)]{LWCDAJa} and \eqref{Text=ext} that
  the complex $\RHom{M}{N}$ has bounded below homology. Moreover, each
  homology module $\H[-i]{\RHom{M}{N}} = \Ext{i}{M}{N}$ is finitely
  generated, as $M$ and $N$ are finitely generated. Now \eqref{wdtinf}
  and \prpref[Propositions~]{wdtI} and \prpref[]{wdtP} yield
  \begin{align*}
    \sup\setof{i \in \ZZ}{\Ext{i}{M}{N}\ne 0} &= -\inf{\RHom{M}{N}}\\
    &= -\wdt{\RHom{M}{N}}\\
    &= -\dpt{M} - \wdt{N} + \dptR\\
    &= \dptR - \dpt{M}.\qedhere
  \end{align*}
\end{bfhpg*}

\section*{Acknowledgments}

\noindent
We thank Olgur Celikbas and Srikanth Iyengar for their comments on an
earlier version of this paper. We also acknowledge an anonymous
referee whose pertinent suggestions helped us tighten the presentaion.

\bibliographystyle{amsplain}

\def\cprime{$'$}
\providecommand{\arxiv}[2][AC]{\mbox{\href{http://arxiv.org/abs/#2}{\sf
      arXiv:#2 [math.#1]}}}
\providecommand{\oldarxiv}[2][AC]{\mbox{\href{http://arxiv.org/abs/math/#2}{\sf
      arXiv:math/#2
      [math.#1]}}}\providecommand{\MR}[1]{\mbox{\href{http://www.ams.org/mathscinet-getitem?mr=#1}{#1}}}
\renewcommand{\MR}[1]{\mbox{\href{http://www.ams.org/mathscinet-getitem?mr=#1}{#1}}}
\providecommand{\bysame}{\leavevmode\hbox to3em{\hrulefill}\thinspace}
\providecommand{\MR}{\relax\ifhmode\unskip\space\fi MR }
\providecommand{\MRhref}[2]{%
  \href{http://www.ams.org/mathscinet-getitem?mr=#1}{#2} }
\providecommand{\href}[2]{#2}

\end{document}